\documentclass[english]{article}
\usepackage[T1]{fontenc}
\usepackage[latin9]{inputenc}
\usepackage{amssymb}

\makeatletter
\newcommand{\lyxaddress}[1]{
\par {\raggedright #1
\vspace{1.4em}
\noindent\par}
}

\date{}

\makeatother

\usepackage{babel}
\begin{document}

\title{From $e$ to $\pi$: \\
Derivation of the Wallis formula for $\pi$ from $e$ }

\author{Ali Sanayei}

\maketitle

\lyxaddress{\begin{center}
\emph{\small{}Center for Optical Quantum Technologies, Institute
for Laser Physics,}\\
\emph{\small{}University of Hamburg, Luruper Chaussee 149, D-22761
Hamburg, Germany}\\
E-mail: \texttt{asanayei@physnet.uni-hamburg.de}~\\

\par\end{center}}
\begin{abstract}
In this Note, we start off with the primary representation of $e$
and from there present an elementary short proof for the Wallis formula
for $\pi$.\\
\\

\end{abstract}
A well-known pre-Newtonian formula for $\pi$, represented as an infinite
product, was derived by John Wallis in 1655 using a method of successive
interpolations \cite{Wallis,Pi source book}:

\begin{equation}
\frac{\pi}{2}=\frac{2\,\,.\,\,2}{1\,\,.\,\,3}\,\,\frac{4\,\,.\,\,4}{3\,\,.\,\,5}\,\,\frac{6\,\,.\,\,6}{5\,\,.\,\,7}\,\,\frac{8\,\,.\,\,8}{7\,\,.\,\,9}\,\cdots\label{Wallis formula 1}
\end{equation}
Different proofs have been offered from the diverse points of view
in the recently past decades. For instance, proofs using the trigonometric
integrals were proposed by Borwein, Borwein, and Stewart \cite{Pi AGM,Stewart}.
Sondow and Weisstein found a derivation by trigonometry \cite{Sondow and Weisstein,Weisstein}.
A proof from the point of view of the geometric means was offered
by W{\"a}stlund \cite{Wastlund}. Miller proposed a probabilistic
proof \cite{Miller}, and Kovalyova derived it from a combinatorial-probabilistic
method \cite{Kovalyov}. More recently, Friedmann and Hagen simply
applied the variational principle of quantum mechanics on the pure
Coulomb potential (hydrogen atom) and showed a very interesting proof
from the standpoint of physics \cite{Friedmann and Hagen}.

Eventuated by the Euler's identity $e^{i\pi}=-1$, that relates $e$
and $\pi$ together, one will be inspired to find a path from $e$
to $\pi$ and vice versa. For instance, given the Wallis formula (\ref{Wallis formula 1}),
Pippenger found an infinite product for $e$ \cite{Pippenger}. 

In this Note, we trace a path from $e$ to $\pi$ by offering an elementary
proof for the Wallis formula (\ref{Wallis formula 1}) only by starting
off with the primary representation of $e$. To do that, we take the
Bernoulli's expression 

\begin{equation}
\lim_{n\rightarrow\infty}\left(1+\frac{1}{n}\right)^{n}=e\label{e}
\end{equation}
as the definition of $e$ (cf. Ref. \cite{Bernoulli}). 

From Eq. (\ref{e}) one can readily write

\begin{equation}
\lim_{n\rightarrow\infty}\,\frac{n}{n+\frac{1}{2}}\left(\frac{n}{n-\frac{1}{2}}\right)^{2n}=e\,,\label{e _2}
\end{equation}
or equivalently

\begin{equation}
\lim_{n\rightarrow\infty}\,\frac{1}{n+\frac{1}{2}}\,\frac{n^{2n+1}}{\left(n-\frac{1}{2}\right)^{2n}}\,\frac{1}{e}=1\,.\label{e _3}
\end{equation}
From Eq. (\ref{e _3}) one can then obtain

\begin{equation}
\lim_{n\rightarrow\infty}\,\frac{1}{n+\frac{1}{2}}\left[\frac{n^{n+\frac{1}{2}}\,e^{-n}}{\left(n-\frac{1}{2}\right)^{n}\,e^{-n+\frac{1}{2}}}\right]^{2}=1\,,\label{e _4}
\end{equation}
which with the help of the Stirling's formula (see Refs. \cite{Stirling}
and \cite{Schopohl} for elementary proofs) and with the property
of gamma functions, $z\,\Gamma\left(z\right)=\Gamma\left(z+1\right)$,
leads to

\begin{equation}
\lim_{n\rightarrow\infty}\,\frac{\Gamma\left(n+1\right)}{\Gamma\left(n+\frac{1}{2}\right)}\,\frac{\Gamma\left(n+1\right)}{\Gamma\left(n+\frac{3}{2}\right)}=1\,.\label{Gamma 1}
\end{equation}
Note that $\Gamma\left(n+\frac{1}{2}\right)$ can be written in terms
of $\Gamma\left(2n\right)$ and $\Gamma\left(n\right)$ by the Legendre
duplication formula, $\Gamma\left(2n\right)=\frac{1}{\Gamma\left(\frac{1}{2}\right)}\,2^{2n-1}\,\Gamma\left(n\right)\Gamma\left(n+\frac{1}{2}\right)$,
see Ref. \cite{Weisstein _ Legendre} for an elementary short proof.
Since $\Gamma\left(\frac{1}{2}\right)=\sqrt{\pi}$ and $\Gamma\left(n\right)=\left(n-1\right)!$
for $n\in\mathbb{N}$, one can then write:

\begin{equation}
\frac{\Gamma\left(n+1\right)}{\Gamma\left(n+\frac{1}{2}\right)}=\frac{1}{\sqrt{\pi}}\,\frac{2^{n}\,n\,\Gamma\left(n\right)}{\frac{2n\,\Gamma\left(2n\right)}{2^{n}\,n\,\Gamma\left(n\right)}}=\frac{1}{\sqrt{\pi}}\,\prod_{j=1}^{n}\,\frac{2j}{2j-1}\,.\label{first part}
\end{equation}
Applying the same argument and writing $\Gamma\left(n+\frac{3}{2}\right)$
in terms of $\Gamma\left(2n+2\right)$ and $\Gamma\left(n\right)$
by the Legendre duplication formula, one obtains:

\begin{equation}
\frac{\Gamma\left(n+1\right)}{\Gamma\left(n+\frac{3}{2}\right)}=\frac{2}{\sqrt{\pi}}\,\frac{2^{n}\,n\,\Gamma\left(n\right)}{\frac{\Gamma\left(2n+2\right)}{2^{n}\,n\,\Gamma\left(n\right)}}=\frac{2}{\sqrt{\pi}}\,\prod_{j=1}^{n}\,\frac{2j}{2j+1}\,.\label{second part}
\end{equation}
Inserting (\ref{first part}) and (\ref{second part}) into (\ref{Gamma 1})
yields

\begin{equation}
\lim_{n\rightarrow\infty}\,\prod_{j=1}^{n}\,\frac{\left(2j\right)\,\left(2j\right)}{\left(2j-1\right)\,\left(2j+1\right)}=\frac{\pi}{2}\,,\label{derivation}
\end{equation}
which is the Wallis formula (\ref{Wallis formula 1}).


\begin{thebibliography}{10}
{\small{}\bibitem{Wallis}J. Wallis, }\emph{\small{}Arithmetica Infinitorium}{\small{}
(Oxford, 1655).}{\small \par}

{\small{}\bibitem{Pi source book}L. Berggern, J. M. Borwein, and
P. B. Borwein, }\emph{\small{}Pi: A Source Book}{\small{} (Springer,
1977).}{\small \par}

{\small{}\bibitem{Pi AGM}J. M. Borwein and P. B. Borwein, }\emph{\small{}Pi
and the AGM}{\small{} (John Wiley \& Sons, 1987).}{\small \par}

{\small{}\bibitem{Stewart}J. Stewart, }\emph{\small{}Calculus: Early
Transcendentals}{\small{} (Cengage Leanring, 2012).}{\small \par}

{\small{}\bibitem{Sondow and Weisstein}J. Sondow and E. W. Weisstein,
``Wallis formula,'' From }\emph{\small{}MathWorld}{\small{}.}\\
{\small{} }\texttt{\small{}http://mathworld.wolfram.com/WallisFormula.html}{\small \par}

{\small{}\bibitem{Weisstein}E. W. Weisstein, ``Pi formulas,'' From
}\emph{\small{}MathWorld}{\small{}.}\\
{\small{} }\texttt{\small{}http://mathworld.wolfram.com/PiFormulas.html}{\small \par}

{\small{}\bibitem{Wastlund}J. W{\"a}stlund, Link{\"o}ping stud.
Math. }\textbf{\small{}2}{\small{}, 1 (2005).}{\small \par}

{\small{}\bibitem{Miller}S. J. Miller, Amer. Math. Monthly }\textbf{\small{}115}{\small{},
740 (2008).}{\small \par}

{\small{}\bibitem{Kovalyov}M. Kovalyov, J. Math. Stat. }\textbf{\small{}5}{\small{},
408 (2009).}{\small \par}

{\small{}\bibitem{Friedmann and Hagen}T. Friedmann and C. R. Hagen,
J. Math. Phys. }\textbf{\small{}56}{\small{}, 112101 (2015).}{\small \par}

{\small{}\bibitem{Pippenger}N. Pippenger, Amer. Math. Monthly }\textbf{\small{}87}{\small{},
391 (1980).}{\small \par}

{\small{}\bibitem{Bernoulli}J. Bernoulli, Journal des S\c{c}avans,
314 (1685).}{\small \par}

{\small{}\bibitem{Stirling}D. Romik, Amer. Math. Monthly }\textbf{\small{}107}{\small{},
556 (2000).}{\small \par}

{\small{}\bibitem{Schopohl}N. Schopohl, }\emph{\small{}Skript zur
Vorlesung Statistische Physik und Thermodynamik}{\small{} (University
of T{\"u}bingen, 2012). }{\small \par}

{\small{}\bibitem{Weisstein _ Legendre}E. W. Weisstein, ``Legendre
duplication formula,'' From }\emph{\small{}MathWorld}{\small{}.}\\
{\small{} }\texttt{\small{}http://mathworld.wolfram.com/LegendreDuplicationFormula.html}{\small \par}\end{thebibliography}
\end{document}